\documentclass[11pt]{amsart}
\usepackage[T1]{fontenc}
\usepackage[english]{babel}
\usepackage{amsmath,amssymb,amsthm,enumitem,xspace}
\usepackage{comment}
\usepackage{accents}
\usepackage[svgnames]{xcolor}
\usepackage[colorlinks=true,citecolor=Indigo,linkcolor=Indigo]{hyperref}
\usepackage{cleveref}
\usepackage{mathrsfs}

\numberwithin{equation}{section}

\newtheorem{thm}{Theorem}

\newtheorem{lem}[thm]{Lemma}

\newtheorem{ithm}{Theorem}

\newtheorem*{ilem*}{Lemma}
\newtheorem*{iconj*}{Conjecture}

\newtheorem*{iprob*}{Problem}

\theoremstyle{definition}

\newtheorem*{defi*}{Definition}
\newtheorem{rem}[thm]{Remark}
\newtheorem*{rem*}{Remark}

\newtheorem*{rems*}{Remarks}
\newtheorem{exam}[thm]{Example}
\theoremstyle{definition}

\newcommand*{\sC}{\mathscr{C}}

\newcommand*{\sS}{\mathscr{S}}

\newcommand*{\ZZ}{\mathbf{Z}}

\newcommand*{\RR}{\mathbf{R}}
\newcommand*{\NN}{\mathbf{N}}

\newcommand*{\se}{\subseteq}

\newcommand*{\lra}{\longrightarrow}
\newcommand*{\Id}{\mathrm{Id}}

\DeclareMathOperator{\Span}{Span}
\newcommand*{\TN}{\langle T\rangle}
\makeatletter
\@namedef{subjclassname@2020}{\textup{2020} Mathematics Subject Classification}
\makeatother

\title[Eigenfunctionals for positive operators]{Eigenfunctionals for positive operators}

\author[Nicolas Monod]{Nicolas Monod}
\address{EPFL, Switzerland}
\email{nicolas.monod@epfl.ch}
\begin{document}

\begin{abstract}
We establish an eigenfunctional theorem for positive operators, evocative of the Krein--Rutman theorem. A more general version gives a joint eigenfunctional for commuting operators.
\end{abstract}
\maketitle
\thispagestyle{empty}


\section{Introduction}
The following ``eigenfunctional theorem'' is an abstract algebraic counterpart to the fundamental Krein--Rutman eigenvector theorem~\cite{Krein-Rutman}, itself an infinite-dimensional version of the Perron--Frobenius principle~\cite{Perron,Frobenius1912}.

\begin{ithm}\label{thm:KR}
Let $V$ be a non-singular ordered vector space and $T\colon V\to V$ any positive linear map.

For every $v\gneqq 0$ with order-bounded orbit $\TN v$ there is a positive linear functional $J$ on the orbit ideal $V_{\TN v}$ and a scalar $t \geq 0$ such that
\begin{equation*}
J(T u) = t\, J(u)  \kern3mm (\forall\, u\in V_{\TN v})  \kern3mm  \text{and} \kern3mm J(v)=1.
\end{equation*}
Moreover, we can choose $J$ such that $t\leq 1$.
\end{ithm}

A first motivation for \Cref{thm:KR} is that it enters the proof of a new characterization of amenable groups 
established in~\cite{Monod_vector-pricings_pre}. That criterion characterizes amenability not by invariance, but by stationarity. Specifically, applying \Cref{thm:KR} to a Markov operator, it is shown that a group $G$ is amenable if and only if $\ell^\infty(G)$ carries a stationary ``conditional mean''; see~\cite{Monod_vector-pricings_pre} for definitions.

\medskip

In the above theorem, $V$ is \textbf{non-singular} if for every $v\gneqq 0$ in $V$ there exists some positive linear functional which is non-zero on $v$. This holds in most familiar examples; for instance, it holds as soon as $V$ admits a Hausdorff locally convex topology compatible with the order~\cite[V.4.1]{Schaefer}. We emphasize however that \Cref{thm:KR} makes no topological assumption; this is a notable difference with the Krein--Rutman theorem and its generalizations, which consider \emph{compact} continuous linear operators on topological vector spaces.

As for the \textbf{orbit ideal}, it refers to the smallest vector subspace containing the orbit $\TN v=\{T^n v : n\in\NN\}$ and closed under order intervals. The orbit is said \textbf{order-bounded} if it admits some upper bound in $V$.

\begin{rem}\label{rem:KR}
The mild assumptions of \Cref{thm:KR} cannot be dropped. First, without the non-singularity assumption there might not exist any non-zero positive functional $J$ at all. Secondly, $J$ cannot in general be defined on all of $V$; consider e.g.\ $v=\delta_0$ (the Kronecker delta) in $V=\ell^\infty(\ZZ)$ for the right shift $T$. Lastly, regarding the boundedness of the orbit, one could guess that it is assumed to ensure $t\leq 1$, but the situation is more interesting than that:

(i)~Even if arbitrary eigenvalues $t$ are allowed, the order-boundedness cannot be dropped completely. See \Cref{exam:no-Krein-Rutman} below.

(ii)~Assuming order-bounded orbits, the additional statement $t\leq 1$ is still not automatic. Consider again the right shift $T$ on $V=\ell^\infty(\ZZ)$ with $v=\delta_0$. The orbit ideal $V_{\TN v}$ consists of the finitely supported functions on $\NN$. Now every $t\geq 0$ admits a functional $J$ as above, namely $J(u) = \sum_{n=0}^\infty t^n\, u(n)$, which is a finite sum.
\end{rem}

\Cref{thm:KR} can also be viewed as a non-compact version of the Markov--Kakutani theorem~\cite{Markov_FP,Kakutani38}. Just like for Markov--Kakutani, we can more generally consider a commuting family of operators instead of a single one. We then obtain a joint eigenfunctional:

\begin{ithm}\label{thm:multi-KR}
Let $V$ be a non-singular ordered vector space and let $S$ be a commutative monoid of positive linear maps $V\to V$.

For every $v\gneqq 0$ with order-bounded $S$-orbit there is a positive linear functional $J$ on the orbit ideal $V_{S v}$ and a monoid homomorphism $t\colon S \to [0,1]$ such that
\begin{equation*}
J(T u) = t(T)\, J(u)  \kern3mm (\forall\, T\in S \ \forall\, u\in V_{S v})  \kern3mm  \text{and} \kern3mm J(v)=1.
\end{equation*}
\end{ithm}

There are however fundamental differences between our setting of positive operators and the convex compact setting of Markov--Kakutani (with linear or affine maps). For fixed points in convex compact spaces, it is straightforward to go from commutativity to, for instance, solvable groups. In fact, a group or monoid satisfies the fixed-point property for convex compact spaces \emph{if and only if} it is amenable. This fails dramatically for positive operators. If we consider the property of admitting an eigenfunctional on the orbit ideal of any order-bounded orbit in a non-singular ordered vector space, already metabelian groups can fail that property, see~\cite{Monod_vector-pricings_pre} and~\cite{Monod_cones} for an extensive discussion.

\emph{We do not even know whether the product of two groups or monoids with that property still satisfies it. We conjecture it does not.}

\smallskip
These stark differences are due to the fact that Theorems~\ref{thm:KR} and~\ref{thm:multi-KR} rely on a combination of positivity and compactness techniques, as opposed to only compactness. This is analogous --- and in fact related --- to the differences between Martin and Poisson boundaries in the theory of Markov processes.

\medskip

Finally, we state a technical variant of \Cref{thm:multi-KR}. It is more cumbersome, but slightly more general and it turns out to have applications to the study of locally compact groups (which we plan to present in a forthcoming paper).

\begin{ithm}\label{thm:multi-KR-bis}
Let $S$ be a commutative monoid of positive linear maps $V\to V$ in an ordered vector space $V$ and pick $v\in V^+$.

Suppose that on the space $\Span(S v)$ spanned by the orbit $S v$ there is a positive linear functional that is bounded on $S v$ and non-zero on $v$.

Then there is a positive linear functional $J$ on the orbit ideal $V_{S v}$ and a monoid homomorphism $t\colon S \to [0,1]$ such that
\begin{equation*}
J(T u) = t(T)\, J(u)  \kern3mm (\forall\, T\in S \ \forall\, u\in V_{S v})  \kern3mm  \text{and} \kern3mm J(v)=1.
\end{equation*}
\end{ithm}

The assumptions of \Cref{thm:multi-KR-bis} are satisfied in the setting of \Cref{thm:multi-KR}. Indeed suppose that the orbit $S v$ of $v\gneqq 0$ is order-bounded by some $w\in V$ and that $V$ is non-singular. Then there is a positive linear functional $\zeta$ on $V$ with $\zeta(v)\neq 0$. Consider the restriction of $\zeta$ to $\Span(S v)$; this is bounded by $\zeta(w)$ on the orbit $S v$.

\medskip
To conclude this introduction, we illustrate that we cannot drop completely from \Cref{thm:KR} the assumption that the orbit is order-bounded, even if we allow arbitrary eigenvalues and keep all other assumptions.

\begin{exam}\label{exam:no-Krein-Rutman}
Let $V$ be the ordered vector space of all finitely supported functions on $\NN$; it is non-singular and its dual can be identified with the space $V^*=\RR^\NN$ of all functions. Consider the positive linear operator $T\colon V\to V$ defined by $T(\delta_n)= \delta_0 + (n+1) \delta_{n+1}$ on the standard basis and take $v=\delta_0$. Then the orbit ideal $V_{\TN v}$ is all of $V$. We claim that regardless of the value $t\in \RR$, possibly greater than one, there is no positive linear functional $J$ on $V$ with $J(v)=1$ and such that $J(T u) = t J(u)$ holds for all $u\in V$.

Suppose otherwise and consider the non-negative sequence defined by $j_n = J(\delta_n)$. Then $j_0=1$ and we have the recurrence relation $j_n=(t j_{n-1} -1)/n$ for all $n\geq 1$. This shows that we have
\begin{equation*}
j_n = \frac1{n!}\left( t^n - \sum_{k=1}^n (k-1)! \,t^{n-k} \right) \kern3mm(\forall\, n\geq 1).
\end{equation*}
Since each $j_n$ is non-negative and $t=J(Tv)\geq 0$, this implies in particular $t^n -(n-1)! \geq 0$ for all $n$. This cannot hold for any $t\in\RR$.
\end{exam}

In view of the formulation of \Cref{thm:multi-KR-bis}, it must in fact be the case that no positive linear functional $\zeta$ on $V$ with $\zeta(\delta_0)\neq 0$ is bounded on the orbit $\{T^n(\delta_0): n\in \NN\}$. This can also be verified directly by showing inductively that the value at~$0$ of $T^n(\delta_0)$ satisfies
\begin{equation*}
\forall n\geq 0: \kern3mm \big(T^{n+1}(\delta_0)\big)(0) = \sum_{p=0}^n (n-p)!\, \big(T^{p}(\delta_0)\big)(0),
\end{equation*}
which defines the super-exponential sequence $1,1,2,5,15,54,235,1237,7790,\ldots$ (the OEIS sequence A051295); we spare the details. It follows that
\begin{equation*}
\zeta\big(T^{n}(\delta_0)\big) \geq \big(T^{n}(\delta_0)\big)(0) \cdot \zeta(\delta_0)
\end{equation*}
is unbounded, indeed super-exponential.  (The latter also gives another proof that no eigenfunctional $J$ exists.)

\section{Proof of the theorems}\label{sec:proofs}
\begin{flushright}
\begin{minipage}[t]{0.85\linewidth}\itshape\small
Obwohl dies ein rein algebraischer Satz ist, so ist mir doch sein Beweis mit den gew{\"o}hnlichen Hilfsmitteln der Algebra nicht gelungen.
\begin{flushright}
\upshape Oskar Perron~\cite{Perron}, p.~262.
\end{flushright}

\smallskip
\upshape
[Although this is a purely algebraic theorem, I have not succeeded in proving it with the usual tools of algebra.]
\end{minipage}
\end{flushright}
We have no topology on the abstract ordered vector space $V$ and thus no topological assumptions on our positive linear maps $V\to V$. However, any (algebraic) dual $V^*$ can be endowed with the weak-* topology, whose uniform structure is complete. This is because $V^*$ is isomorphic, as a topological vector space, to the product of copies of $\RR$ indexed by a basis of $V$, see e.g.~\cite[II \S\,6 No.\,6--7]{BourbakiEVT81}.

We will also consider this for the ideal $V_{A}$ instead of $V$, where $A=\TN v$ and more generally $A=S v$. We recall the straightforward observation that for any set $A\se V^+$ the ideal $V_A$ can be written explicitly as
\begin{equation*}\tag{$*$}
V_A = \big\{ u\in V : \exists a_1, \ldots, a_q\in A \text{ with } \pm u \leq a_1 + \cdots + a_q \big\}.
\end{equation*}
It follows that the order of $V$ restricted to $V_A$ is \textbf{generating}, i.e.\ $V_A=\Span(V_A^+)$. In particular, the dual positive cone $(V_A^*)^+$ in the algebraic dual $V_A^*$ is \textbf{proper}, which means  $-(V_A^*)^+\cap (V_A^*)^+ = \{0\}$.

The proof of the eigenfunctional theorems will involve extremal rays. Recall that even very familiar cones may fail to contain any extremal ray: for instance, the positive cone in the Hilbert space $L^2[0,1]$, which is weakly complete and self-dual, has no extremal ray.  Choquet's fundamental insight is that a cone $C$ in a locally convex topological vector space will be the closed convex hull of its extremal rays as soon as the following holds: 

\itshape
for every $\lambda \in C$ there is a compact convex set $B\se C$ containing~$\lambda$ and such that $C\smallsetminus B$ is convex too.\upshape

Such a set $B$ is what Choquet calls a \textbf{chapeau} (a hat, or cap). Indeed, it then suffices to apply the Krein--Milman theorem to every chapeau. We refer to~\cite{Choquet63} for an overview and to~\cite{Choquet_LAII} or~\cite[II \S\,7]{BourbakiEVT81} for the full monty.

When the above condition holds, the cone is called \textbf{well-capped} (bien coiff{\'e}). It is apparent on the definition that this property is inherited by closed subcones.

We shall need the following (rather standard) fact, in which the countability assumption is essential.

\begin{lem}\label{lem:chapeau}
Let $V$ be an ordered vector space and $A\se V^+$ a countable set. Let $C= (V_A^*)^+$ be the cone of positive linear functionals in the algebraic dual $V_A^*$ of the ideal generated by $A$, endowed with the weak-* topology. 

Then $C$ is well-capped.
\end{lem}

\noindent
Without countability, a counter-example is $V=L^2[0,1]$ with $A=V^+$. In that case, $V_A=V$ and $C=(V_A^*)^+\cong V^+$ (e.g.\ by Klee's theorem~\cite[Thm.~B]{Klee55b}). Thus $C$ has no extremal rays and hence is not well-capped for any topology.

\begin{proof}[Proof of \Cref{lem:chapeau}]
We fix $\lambda\in C$, which we can assume non-zero. We also choose an enumeration $A=\{v_n : n\in \NN\}$. (If $A$ is finite there will be repetitions, but that case is trivial anyways, since then $C$ admits a compact base). We claim that there is $n$ with $\lambda(v_n)\neq 0$. Indeed, choose $u\in V_A$ with $\lambda(u)\neq 0$. There are indices $n_1, \ldots, n_q$ with $\pm u \leq v_{n_1} + \cdots + v_{n_q}$; the claim follows by positivity of $\lambda$.

This implies that we can choose a sequence of scalars $t_n>0$ such that $\sum_{n\geq 0} t_n \lambda(v_n) = 1$. We define a map
\begin{equation*}
s\colon C \lra [0, +\infty], \kern3mm s(\xi) = \sum_{n\geq 0} t_n \xi(v_n).
\end{equation*}
This map is additive and positively homogeneous on $C$, where this is suitably understood for $[0, +\infty]$-valued maps. We now consider the sublevel set $B=\{\xi\in C : s(\xi) \leq 1\}$. The only millinerial property of $B$ that is not clear from this definition is compactness.

Since $B$ is contained in the complete space $C$, it suffices to justify that it is closed and bounded, both with respect to the weak-* topology~\cite[23.11]{Choquet_LAII}. We recall here that a set is \emph{bounded} in the sense of topological vector spaces if it is absorbed by a sufficiently large rescaling of any neighbourhood of~$0$. Every partial sum $\xi\mapsto  \sum_{n= 0}^N t_n \xi(v_n)$ is weak-* continuous, and therefore the map $s$ is lower semi-continuous; this shows that $B$ is closed.

For boundedness, after unwinding the definitions, it suffices to show that for every $u_1, \ldots, u_m\in V_A$ there is $r>0$ such that
\begin{equation*}
B \se \big\{\xi\in C : |\xi(u_j)| < r \ \forall\, j \big\}.
\end{equation*}
We thus fix some $u_1, \ldots, u_m\in V_A$ and $\xi\in B$. Using~$(*)$ for each $u_j$ and grouping all sums together, we obtain indices $n_1, \ldots, n_q$ with
\begin{equation*}
\forall\, j: \kern3mm \pm u_j \leq v_{n_1} + \cdots + v_{n_q}.
\end{equation*}
Given $i=1, \ldots , q$ we have
\begin{equation*}
\xi(v_{n_i}) = \frac1{t_{n_i}} t_{n_i}\xi(v_{n_i}) \leq \frac1{t_{n_i}} s(\xi) \leq \frac1{t_{n_i}},
\end{equation*}
which implies
\begin{equation*}
 |\xi(u_j)| \leq \xi(v_{n_1}) + \cdots + \xi(v_{n_q}) \leq  \sum_{i=1}^q  \frac1{t_{n_i}}.
\end{equation*}
Thus we can choose any $r>  \sum_{i=1}^q  \frac1{t_{n_i}}$ to conclude the proof.
\end{proof}

We now prove the theorems stated in the introduction. We have already recorded that \Cref{thm:multi-KR-bis} implies \Cref{thm:multi-KR} which in turn contains \Cref{thm:KR}.

Since \Cref{lem:chapeau} fails without the countability assumption, it is natural to prove \Cref{thm:multi-KR-bis} first for countable monoids $S$. This already contains the case of \Cref{thm:KR}.

\begin{proof}[Proof of \Cref{thm:multi-KR-bis} for $S$ countable]
We recall the setting. We are given a countable commutative monoid $S$ of positive linear maps $T\colon V\to V$ on an ordered vector space $V$, an element $v\geq 0$ and a positive linear functional $\zeta$ on the orbit span $\Span(S v)$ such that $\zeta$ is bounded on $S v$ and $\zeta(v)\neq 0$. In particular, $v\gneqq 0$.

We first note that $\zeta$ can be extended to a positive linear functional on the orbit ideal $V_{S v}$ by the Kantorovich theorem (\cite[Theorem~1.6.1]{Jameson} or~\cite[Theorem~1.36]{Aliprantis-Tourky}). That theorem requires $\Span(S v)$ to majorize $V_{S v}$ (i.e. to be cofinal in it), which is the case by~$(*)$.

Each $T$ restricts to a positive linear map $V_{S v}\to V_{S v}$ on the ideal $V_{S v}$; we still denote it by $T$. We endow the algebraic dual $V_{S v}^*$ of $V_{S v}$ with the weak-* topology, which is complete as mentioned earlier. It follows that the positive cone $C=(V_{S v}^*)^+$ is also complete and hence weakly complete since the weak-* topology is a weak topology.

Consider the monoid $S^*$ of all dual maps $T^* \colon V_{S v}^* \to V_{S v}^*$ as $T$ ranges over $S$. That is, $(T^* \xi)(u) = \xi(T u)$ for $u\in V_{S v}$ and $\xi\in V_{S v}^*$. Algebraically, this is the opposite monoid to (the restriction to $V_{S v}$ of ) $S$ and in particular it is still a commutative monoid. Each $T^*$ is a positive linear map and $S^*$ leaves $C$ invariant, i.e.\ $T^*(C) \se C$ for all $T\in S$.

\medskip
We claim that $C$ contains a closed $S^*$-invariant subcone $C_1\se C$ which is minimal among all closed $S^*$-invariant subcones satisfying the following property:
\begin{center}\itshape
$(\dagger)$ there is $\xi\in C_1$ with $\xi(v)=1$ and with $\xi|_{S v} \leq 2$.
\end{center}
We first show that any closed $S^*$-invariant subcone $C_0\se C$ will satisfy~$(\dagger)$ as soon as it satisfies the formally weaker property
\begin{center}\itshape
$(\ddagger)$ there is $\xi\in C_0$ with $\xi(v)=1$ and which is bounded on $S v$.
\end{center}
Indeed, $(\ddagger)$ allows us to define $s=\sup_{T \in S} \xi(T v) \geq 1$ and to choose $T_0\in S$ with $\xi(T_0 v) \geq s /2$. Next, we observe
\begin{equation*}
\sup_{T\in S} (T_0^* \xi) (T v) = \sup_{T\in S} \xi (T_0  T v) \leq \sup_{T\in S} \xi (T v) =s.
\end{equation*}
Therefore, since $C_0$ is $S^*$-invariant, we can replace $\xi$ by
\begin{equation*}
\frac1{(T_0^* \xi)(v)} T_0^* \xi = \frac1{\xi(T_0 v) } T_0^* \xi
\end{equation*}
which satisfies~$(\dagger)$.

Next, we record that $C$ itself enjoys~$(\ddagger)$ by assumption (upon rescaling $\zeta$) and hence also~$(\dagger)$.

Now, to establish the claim, it suffices by Zorn's lemma to show that for any chain $\sC$ of closed $S^*$-invariant subcones $K\se C$ containing an element $\xi_K$ as in~$(\dagger)$, the intersection $\bigcap_{K\in \sC} K$ still contains such an element.

For every $u\in V_{S v}$, we can by~$(*)$ fix $q\in \NN$ and $T_1, \ldots, T_q\in S$ such that $\pm u \leq T_1 v + \cdots + T_q v$. This gives for any $K\in\sC$ the bound
\begin{equation*}
  |\xi_K (u)| \leq \xi_K ( T_1 v) + \cdots + \xi_K ( T_q v) \leq 2 q,
\end{equation*}
which is uniform over $K\in\sC$. This implies, by the criterion given in~\cite[23.11]{Choquet_LAII}, that the net $(\xi_K)_{K \in \sC}$ is weak-* relatively compact in $C$, recalling that $C$ is weak-* complete. We can therefore choose an accumulation point $\xi\in C$ of this net.

We have $\xi\in \bigcap_{K\in\sC} K$ since each $K$ is closed. Moreover, the properties $\xi_K(v)=1$ and $\xi_K (T v) \leq 2$ (for any $T\in S$) pass to weak-* limits. This completes the proof of the claim.

\medskip
Let thus $C_1$ be minimal as claimed and choose some $T\in S$. The set
\begin{equation*}
C_T = (\Id + T^*) C_1  =  \{\eta + T^* \eta : \eta \in C_1\}
\end{equation*}
is a  subcone of $C_1$. It is $S^*$-invariant because $S^*$ is commutative; indeed, for any $ T'\in S$ we have
\begin{equation*}
T'^*C_T = (T'^* +T'^* T^*) C_1 = (\Id + T^*) (T'^*C_1) \se  (\Id + T^*) C_1 = C_T.
\end{equation*}
Let $\xi\in C_1$ be a witness for~$(\dagger)$. Then $\xi + T^* \xi \in C_T$ is bounded on $S v$ and non-zero on $v$; therefore, after rescaling, it witnesses that $C_T$ satisfies~$(\ddagger)$ and hence also~$(\dagger)$. 

We claim that $C_T$ is closed. Indeed, if any net $(\eta_\alpha + T^* \eta_\alpha)$ (where $\alpha$ ranges over some directed index set) with $\eta_\alpha\in C_1$ converges in $C$, then so does the net $(\eta_\alpha)$ upon passing to a subnet because the sum map $C_1\times C_1 \to C_1$ is proper, see Cor.~1 in~\cite[II \S\,6 No.\,8]{BourbakiEVT81}. (For convergence in terms of properness, see Thm.~1(c) in~\cite[I \S\,10 No.\,2]{BourbakiTGI_alt}.) This implies that the limit of $(\eta_\alpha + T^* \eta_\alpha)$ belongs to $C_T$.

In conclusion, the minimality of $C_1$ shows that we have $C_T=C_1$, where $T\in S$ is arbitrary.

\medskip
By \Cref{lem:chapeau}, the cone $C$ is well-capped and therefore so is $C_1$. It follows that $C_1$ is the closed convex hull of its extremal rays~\cite[II \S\,7 No.\,2 Cor.~2]{BourbakiEVT81}. In particular, extremal rays exist. Consider any $0 \neq J\in C_1$ such that the ray $\RR^+ J$ is extremal. For each $T\in S$, the equality $C_T=C_1$ shows that there is $\eta\in C_1$ with $J= \eta + T^* \eta$. By extremality, both summands are on the ray $\RR^+ J$ and thus we have $a,b\geq 0$ with $\eta=a J$ and $T^* \eta = b J$. The relation $J= \eta + T^* \eta$ implies that $a\neq 0$ and that $b = 1-a$. It follows $T^* J = \frac{1-a}a J$ and we set $t_J(T):= (1-a)/a \geq 0$.

Thus, to each $J$ spanning an extremal ray of $C_1$ corresponds a map $t_J\colon S \to \RR_{\geq 0}$ with
\begin{equation*}
\forall\, T\in S \ \forall\, u\in V_{S v} : \kern3mm J(T u) = t_J(T)\, J(u).
\end{equation*}
We have $J(v)\neq 0$ since otherwise the above property would force $J$ to vanish, being a positive functional on the ideal $V_{S v}$. After rescaling, we assume $J(v)=1$. Next, note that $t_J$ is a monoid homomorphism to the multiplicative monoid $\RR_{\geq 0}$; indeed,
\begin{multline*}
t_J(T_1)  t_J(T_2) =  t_J(T_1)  t_J(T_2) \, J(v) = t_J(T_1)\, J(T_2 v) = J(T_1 T_2 v) \\ =  t_J(T_1 T_2)\, J(v) =t_J(T_1 T_2).
\end{multline*}
It remains to show that we can choose $J$ such that $t=t_J\leq 1$; as we noted in \Cref{rem:KR}, this is not the case for all eigenfunctionals.

\medskip
To that end, we enumerate $S=\{T_i : i\in \NN\}$ and our final claim is that for every $n\in \NN$ there exists $J_n$ as above and such that $t_{J_n}(T_i) \leq 2$ for all $i=1 , \ldots , n$. Supposing this fails for some $n$, we pick $p\in\NN$ such that $2^p\geq 2 n+1$. Then for each $J$ as above there is $i\leq n$ with $t_J(T_i^p) \geq 2 n+1$. Since $t_J\geq 0$, this implies 
\begin{equation*}
J\big(T_1^p v + \cdots + T_n^p v \big) \geq (2 n+1) \, J(v)
\end{equation*}
and thus every $ J\in C_1$ on any extremal ray of $C_1$ belongs to the closed convex set
\begin{equation*}
\Big\{ \zeta \in C_1 : \zeta\big(T_1^p v + \cdots + T_n^p v \big) \geq (2 n+1) \, \zeta(v) \Big\}.
\end{equation*}
On the other hand, $\xi$ does not belong to that set since $\xi(T_j^p v) \leq 2 \xi(v) = 2$ for all $j$. Let $B$ be a chapeau of $C_1$ containing $\xi$ and recall that every extremal point of $B$ lies on an extremal ray of $C_1$ (this is the motivating definition of a chapeau~\cite[Lem.~2]{Choquet62b}). This contradicts the Krein--Milman theorem, confirming the claim.

\medskip
From the final claim, we can deduce that there is $J$ as above but with  $t_J(T) \leq 2$ for all $T\in S$, obtained as any accumulation point of the sequence $(J_n)$. Indeed, by the same compactness principle as for the initial claim, it suffices to justify that for every given $u\in V_{S v}$ the values $J_n(u)$ remain bounded. We write again $\pm u \leq T_{i_1} v + \cdots +  T_{i_q} v$ and deduce that $|J_n(u)| \leq  2 q$ holds for all $n$ large enough, as desired.

Now, $t:= t_J$ being a monoid homomorphism, $t\leq 2$ forces in fact $t(T) \leq 1$ by considering $t(T^n)$ as $n\to\infty$. This concludes the proof of \Cref{thm:multi-KR-bis} in the countable case.
\end{proof}

It remains to reduce the general version of \Cref{thm:multi-KR-bis} to the case of countable monoids. The natural reduction strategy is to consider all countable submonoids $S_0\se S$ and to extend functionals on $V_{S_0} v$ to the larger space $V_{S v}$ in a coherent way. Contrary to the context of convex compact sets, such extensions are not always possible in the setting of positive linear maps; for instance, Hahn--Banach-like extension principles such as the Kantorovich theorem do not hold without further assumptions.

Nevertheless, a compactness argument will get around this difficulty, as follows.

\begin{proof}[Proof of \Cref{thm:multi-KR-bis}, general case]
Given any submonoid $S_0\se S$, we consider the following set $E(S_0)\se  V_{S_0 v}^*$ of positive eigenfunctionals on $V_{S_0 v}$:
\begin{equation*}
E(S_0) = \left\{  J\in (V_{S_0 v}^*)^+:
\begin{array}{l}
J(v)=1 \text{ and }\forall\, T\in S_0, \  J(T v) \leq 1  \text{ and }\\
\forall\, T\in S_0\, \forall\, u\in V_{S_0 v}, \  J(T u) = J(T v) \, J(u) 
\end{array}
\right\}.
\end{equation*}
The countable case of \Cref{thm:multi-KR-bis} states that this set is non-empty for all countable $S_0$ since, in the notation of the theorem, $t(T) = J(T v)$. We need to prove that $E(S)$ is non-empty, too.

We denote by $\sS$ the collection of all countable submonoids $S_0\se S$. Viewing $\sS$ as a directed set for the inclusion, $S$ is the directed union (=colimit) of all $S_0\in \sS$. Moreover, $V_{S v}$ is the directed union of all $V_{S_0 v}$ because of the characterisation~$(*)$ of the orbit ideal $V_{S v}$. Therefore, the algebraic dual  $V_{S v}^*$ is the inverse limit of the algebraic duals $V_{S_0 v}^*$ along $S_0 \in \sS$.

If $S_1$ contains $S_0$, then $ V_{S_1 v}$ contains $V_{S_0 v}$ and the restriction map $V_{S_1 v}^* \to V_{S_0 v}^*$ maps $E(S_1)$ to $E(S_0)$ (the fact that it is a priori not surjective is the extension failure mentioned before the proof). Thus we can form the inverse limit $E=\lim_{S_0\in \sS} E(S_0)$ and $E$ lies in $V_{S v}^*$. Moreover, $E$ is none other than $E(S)$ because the definition of $E(S)$ is purely in terms of its values on various elements of $V_{S v}$.

In conclusion, it remains only to justify that the inverse limit $E$ is non-empty. The set $E(S_0)$ is compact in the weak-* topology of $V_{S_0 v}^*$ because it is closed and we can bound $|J(u)|$ by $q$ whenever $\pm u \leq T_1 v + \cdots + T_q v$.
Moreover, the restriction map $V_{S_1 v}^* \to V_{S_0 v}^*$ is dual and hence weak-* continuous. In conclusion, $E$ is the inverse limit of a system of non-empty compact Hausdorff spaces, which implies that $E$ is non-empty.
\end{proof}

\begin{rem}
Some of the ideas in the proofs above are related to our earlier work on fixed points for group actions on cones~\cite{Monod_cones}. What makes the present situation very different is that a general positive map, typically non-invertible, will often not fix any non-zero point.

This difference is crucially important for the applications to amenability considered in~\cite{Monod_vector-pricings_pre}, where the Markov operator will not fix a non-zero functional unless the group is amenable.
\end{rem}


\bibliographystyle{amsalpha-nobysame}
\bibliography{../BIB/ma_bib}

\end{document}